\documentclass[11pt, a4paper, leqno]{article}
\usepackage[intlimits]{amsmath}
\usepackage{amsthm,amssymb}
\usepackage{amsfonts}
\usepackage{fancyhdr}
\usepackage{verbatim}
\usepackage{color}
\definecolor{labelkey}{rgb}{1,0,0}
\usepackage{graphicx}
\DeclareMathOperator{\im}{Im}
\DeclareMathOperator{\re}{Re}

\newcommand{\Pcal}{\mathcal{P}}
\newcommand{\Scal}{\mathcal{S}}
\newcommand{\Rbb}[1][]{\mathbb{R}^{#1}}
\newcommand{\Cbb}[1][]{\mathbb{C}^{#1}}
\newcommand{\Sbb}{{\mathbb{S}}}
\newcommand{\D}{{\mathbb D}}

\allowdisplaybreaks[4]
\numberwithin{equation}{subsection}
\theoremstyle{plain}
\newtheorem{thm}[equation]{Theorem}

\newtheorem{prop}[equation]{Proposition}

\theoremstyle{definition}

\newtheorem{conj}[equation]{Conjecture}

\theoremstyle{remark}

\newcommand{\pd}[2]{\frac{\partial{#1}}{\partial{#2}}}

\newcommand{\ii}{\text{i}}
\newcommand{\id}{\text{Id}}
\renewcommand{\tilde}{\widetilde}
\newcommand{\DD}{(\mathcal{D}\mathcal{D})_\lambda}
\newcommand{\DN}{(\mathcal{D}\mathcal{N})_\lambda}
\newcommand{\ND}{(\mathcal{N}\mathcal{D})_\lambda}
\newcommand{\NN}{(\mathcal{N}\mathcal{N})_\lambda}

\newcommand{\cP}{{\Pcal}}

\newcommand{\area}{{\rm Area}}
\newcommand{\sI}{{\sigma_{\rm{I}}}}
\newcommand{\sII}{{\sigma_{\rm{II}}}}

\begin{document}

\title{Spectral problems with mixed Dirichlet-Neumann boundary conditions: isospectrality and beyond\thanks{The research of M.L. was
partially supported by CRM (Montr\'eal) and by the EPSRC Spectral Theory Network (U.K.).
The research of D.J. was partially supported by NSERC, FQRNT, Sloan Foundation and Dawson Fellowship.
The research of N.N. was partially supported by NSF. The research of I.P. was partially supported
by NSERC, FQRNT and the Edinburgh Mathematical Society.
}}
\author{
Dmitry Jakobson\\
\normalsize\small Dept. of Mathematics and Statistics\\
\normalsize\small McGill University\\
\normalsize\small 805 Sherbrooke Str. West, Montreal\\
\normalsize\small Quebec H3A 2K6 Canada \\
\normalsize\small {\sffamily jakobson@math.mcgill.ca}
\and
Michael Levitin\\
\normalsize\small Department of Mathematics\\
\normalsize\small Heriot-Watt University\\
\normalsize\small Riccarton, Edinburgh EH14 4AS\\
\normalsize\small United Kingdom\\
\normalsize\small {\sffamily M.Levitin@ma.hw.ac.uk}
\and
Nikolai Nadirashvili\\
\normalsize\small Department of Mathematics\\
\normalsize\small University of Chicago\\
\normalsize\small Chicago, IL 60637\\
\normalsize\small USA\\
\normalsize\small {\sffamily nicholas@math.uchicago.edu}
\and
Iosif Polterovich\\
\normalsize\small  D\'ept. de Math\'ematiques et de Statistique\\
\normalsize\small  Universit\'e de Montr\'eal\\
\normalsize\small CP 6128 succ Centre-Ville\\
\normalsize\small  Montreal, QC H3C 3J7, Canada \\
\normalsize\small {\sffamily iossif@dms.umontreal.ca}
}
%
%
%
\date{\small \today}

\maketitle
\begin{abstract} \noindent Consider a bounded domain with the
Dirichlet condition on a
part of the boundary and the Neumann condition on its complement. 
Does the spectrum of the Laplacian determine uniquely which condition 
is imposed on which part? We present some results, conjectures and
problems related to this variation on the isospectral theme.
\end{abstract}
{\small \textbf{Keywords:} Laplacian, mixed Dirichlet-Neumann
problem, isospectrality.}

\

\noindent{\small \textbf{2000 Mathematics Subject Classification:}
58J53, 58J50, 35P15.}
\newpage
\section{Introduction}\label{sec:intro}

Let $\Omega\subset\mathbb{R}^2$ be a bounded domain, its boundary
being decomposed as $\partial\Omega =
\overline{\partial_1\Omega\cup\partial_2\Omega}$, where
$\partial_1 \Omega$, $\partial_2 \Omega$ are finite unions of open
segments of $\partial \Omega$ and $\partial_1 \Omega \cap
\partial_2 \Omega = \emptyset$. Suppose that there are no
isometries of ${\Bbb R}^2$ exchanging  $\partial_1\Omega$ and
$\partial_2\Omega$. We call such a decomposition of the boundary
{\it nontrivial}. Consider a Laplace operator on $\Omega$ and
assume that on one part of the boundary we have the Dirichlet
condition and on the other part the Neumann condition. Does the
spectrum of the Laplacian determine uniquely which condition is
imposed on which part?

We recall the classical question of Mark Kac, ``\emph {Can one
hear the shape of a drum?}" \cite{K} related to the  (Dirichlet)
Laplacian on the plane, which still remains open for smooth (as
well as for convex) domains. For arbitrary planar domains it was
answered negatively in \cite{GWW} using an algebraic construction
of \cite{Sum}, see also reviews and extensions \cite{BCDS, Bro,
Bus} and references therein.

We may reformulate our question in a similar way. Consider two drums with
drumheads which are partially attached to them. The drumhead of the
first drum is attached exactly where the drumhead of the second drum
is free and vice versa. Can one distinguish between the two drums by
hearing them?

Similarly to the question of Kac, the answer to our question is in general 
negative. In this note we construct a family of domains, each of
them having a nontrivial isospectral (with respect to exchanging
the Neumann and Dirichlet boundary conditions) boundary
decomposition. We say that such domains {\it admit
Dirichlet-Neumann isospectrality}.

Our main example is a half-disk which is considered in sections~\ref{subs:ME}--\ref{subs:proof3}. 
In sections~\ref{subs:qs} and~\ref{subs:sb} we construct generalizations of
the main example, in particular to non-planar domains. Section~\ref{subs:nes}
provides a  simple necessary condition for a boundary
decomposition to be Dirichlet-Neumann isospectral. In section~\ref{subs:are}
we discuss if there exist domains not admitting Dirichlet-Neumann
isospectrality and conjecture that a disk should be one of them.
Some numerical evidence in favour of this conjecture is also presented.

Our motivation for the study of Dirichlet-Neumann isospectrality,
rather surprisingly, comes from a seemingly unrelated problem of
obtaining a   sharp upper bound for the first eigenvalue on a
surface of genus two. We discuss it, as well as some other
relevant eigenvalue inequalities,  in sections~\ref{subs:g2}--\ref{subs:bounds}.

\section{Domains isospectral with respect to the Dirichlet-Neumann swap}\label{sec:DN}
\subsection{Main Example}\label{subs:ME}

Our principle example is constructed using the half-disk. Let
$\Omega:=\{z\in\Cbb: |z|<1, \im z>0\}$ be an upper half of a disc
centered at the point $O$ (here and further  on we shall often
identify the real plane $\Rbb[2]$ with the complex plane $\Cbb$
and shall use the complex variable $z=x+\ii y$ instead of real
coordinates $(x,y)$ on the plane). Consider the following boundary
decomposition:
\begin{gather}
\partial_1 \Omega =  \{\re z\in(-1,0),\, \im z=0\}\cup\{|z|=1,\
|\arg z -\pi/2|<\pi/4\}\,,\\
\partial_2 \Omega =  \{\re z\in(0,1),\, \im z=0\}\cup\{|z|=1,\
\pi/4<|\arg z -\pi/2|<\pi/2\}\,.
\end{gather}
Obviously, such a decomposition is non-trivial.

\begin{figure}[thb!]
\begin{center}
\fbox{\resizebox{0.9\textwidth}{!}{\includegraphics*{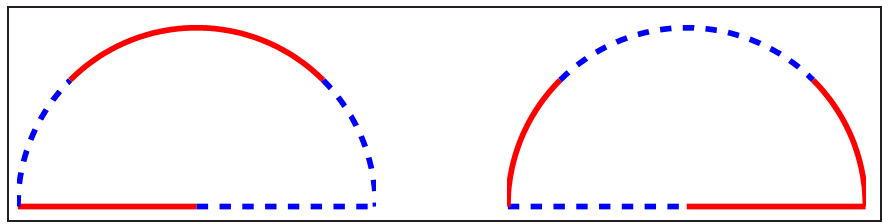}}}
\caption{Problems I and II on the half-disk. Here and further on,
\textcolor{red}{red solid line} denotes the Dirichlet boundary
conditions and \textcolor{blue}{blue dashed line} the Neumann
ones.\label{fig:fig1}}
\end{center}
\end{figure}

Consider the following boundary value spectral problems on $\Omega$:
$$
\text{\textbf{Problem I: }} -\Delta u =\lambda u\quad\text{in
}\Omega\,, \quad u|_{\partial_1 \Omega} =  0\,,\quad
\left.\pd{u}{n}\right|_{\partial_2 \Omega} = 0\,,
$$
and
$$
\text{\textbf{Problem II: }} -\Delta v =\lambda v\quad\text{in
}\Omega\,, \quad v|_{\partial_2 \Omega} = 0\,,\quad
\left.\pd{v}{n}\right|_{\partial_1 \Omega} = 0\,;
$$
(see Figure~\ref{fig:fig1}). Hereinafter $\pd{}{n}$ is the  normal
derivative and $\Delta=\partial^2/\partial x^2 + \partial^2/\partial
y^2$ is the Laplace operator.

Let $\sI$ denote the spectrum of Problem I and $\sII$ --- the
spectrum of Problem II. Both spectra are discrete and positive.
Our main claim is the following

\begin{thm}\label{thm:main}
With account of multiplicities, $\sI\equiv\sII$.
\end{thm}

We give three different proofs of Theorem~\ref{thm:main}, each of
them,  in our opinion, instructive in its own right, and
generalize this Theorem later on for a  certain class of examples.

\subsection{Proof of Theorem~\ref{thm:main} by transplantation}\label{subs:proof1}
This proof uses the transplantation trick similar to \cite{Ber,
BCDS}.  Let $u(z)=u(r,\phi)$ be an eigenfunction of Problem I
corresponding to an eigenvalue $\lambda$; here $(r,\phi)$ denote
the usual polar coordinates. Let us introduce a mapping  $T:
u\mapsto v$, where
\begin{align*}
v(z)&=(Tu)(z)\\
:&=\frac{1}{\sqrt{2}}
\begin{cases}
\label{T}
u(r,\frac{\pi}{2}-\phi)-u(r,\frac{\pi}{2}+\phi)\qquad&\text{if }\phi=\arg z\in(0,\frac{\pi}{2}]\,,\\
u(r,\frac{3\pi}{2}-\phi)+u(r,\phi-\frac{\pi}{2})\qquad&\text{if }\phi=\arg z\in[\frac{\pi}{2},\pi)\,.
\end{cases}
\end{align*}
Then it is easily checked that $v(z)$ is an eigenfunction of
Problem  II: it satisfies the equation and the boundary conditions
as well as the matching conditions for the trace of the function
and the trace of the normal derivative on the central symmetry
line $r\in(0,1)$, $\phi=\pi/2$.

Similarly, if $v(z)$ is an eigenfunction of Problem II, in order
to  construct an eigenfunction $u(z)$ of Problem I we use an
inverse mapping $T^{-1}$ (one may check that $T^8=\id$ and hence
$T^{-1}=T^7$):
\begin{align*}
u(z)&=(T^{-1}v)(z)\\
:&=\frac{1}{\sqrt{2}}
\begin{cases}
v(r,\frac{\pi}{2}-\phi)+v(r,\frac{\pi}{2}+\phi)\qquad&\text{if }\phi=\arg z\in(0,\frac{\pi}{2}]\,,\\
v(r,\frac{3\pi}{2}-\phi)-v(r,\phi-\frac{\pi}{2})\qquad&\text{if }\phi=\arg z\in[\frac{\pi}{2},\pi)\,.
\end{cases}
\end{align*}
This proves that the sets $\sigma_I$ and $\sigma_2$ coincide. The
equality of multiplicities for each eigenvalue follows immediately
from the linearity of the map $T$. \qed

\begin{figure}[thb!]
\begin{center}
\fbox{\parbox{0.2\textwidth}{\resizebox{0.18\textwidth}{!}
{\includegraphics*{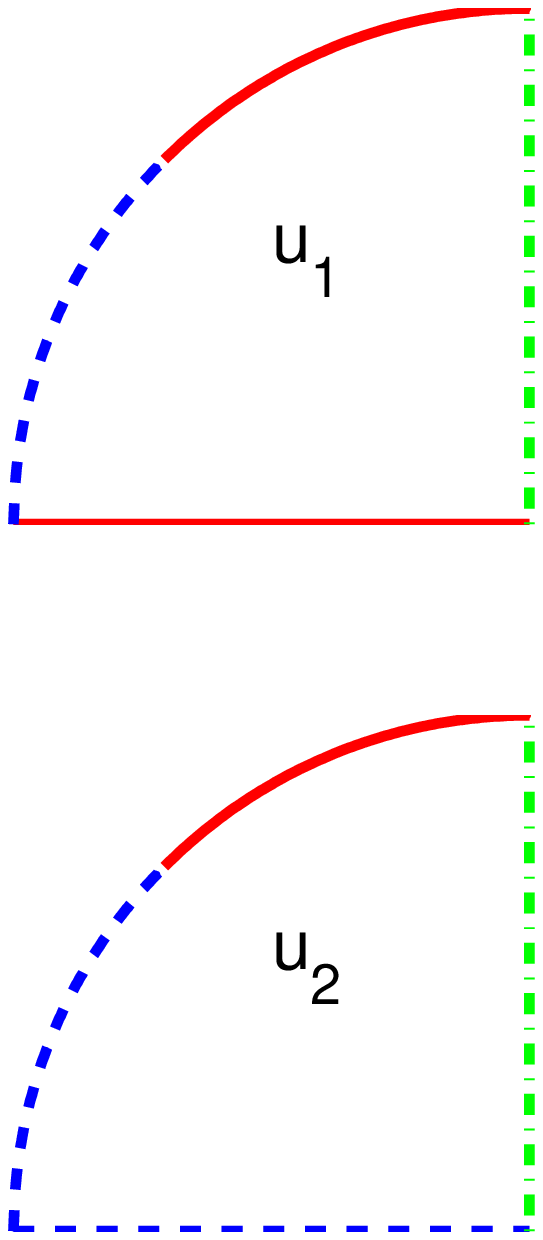}}}
\quad
${u_1\choose u_2}\stackrel{\tilde{T}}
{\longrightarrow}{v_1\choose v_2}=\tilde{T}{u_1\choose u_2}:=\frac{1}{\sqrt{2}}{
u_1-u_2\choose u_1+u_2}$\quad
\parbox{0.2\textwidth}{\resizebox{0.18\textwidth}{!}{\includegraphics*{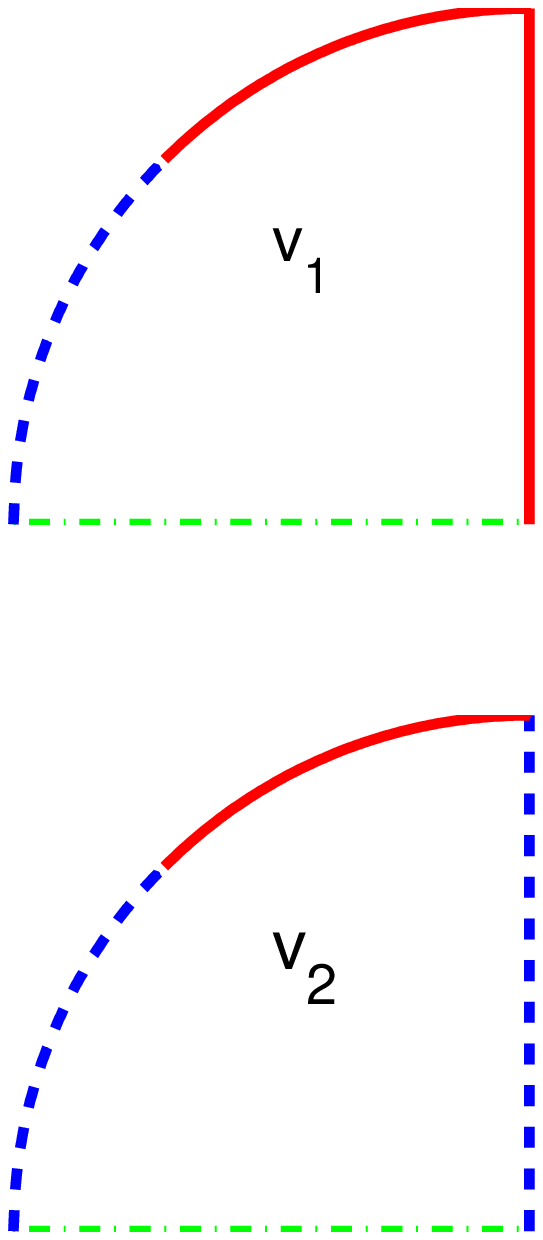}}}}
\caption{Problems $\tilde{\text{I}}$  and $\tilde{\text{II}}$  on the
quarter-disk
and the map $\tilde{T}$.
Here and further on, \textcolor{green}{green dash-dotted line}
denotes the matching conditions.\label{fig:fig2}}
\end{center}
\end{figure}

It is easy to visualise the mapping $T$ in the following way.
Cutting a half-disk along the symmetry line, we can rewrite
Problem I as a system of two  boundary value problems with respect
to functions $u_1$, $u_2$ on a quarter-disk $\Upsilon:=
\{z\in\Cbb: |z|<1\,\ \im z>0\,,\re z<0\}$; $u_1$, $u_2$ should
satisfy the matching conditions ($u_1=u_2$,
$\pd{u_1}{n}=-\pd{u_2}{n}$) on the line $\arg z = \pi/2$, see
Figure~\ref{fig:fig2}. We call this equivalent statement
\textbf{Problem }$\tilde{\text{\textbf{I}}}$. The map $\tilde{T}$
shown in Figure~\ref{fig:fig2} maps the eigenfunction $(u_1,u_2)$
of Problem $\tilde{\text{I}}$ into an eigenfunction $(v_1,v_2)$ of
\textbf{Problem }$\tilde{\text{II}}$, which is in turn equivalent
to Problem II (it is obtained by rotating a half-disk in problem
II clockwise by $\pi/2$ and then  cutting along the symmetry
line).

The possibility of transplanting the eigenfunctions indicates that our
problem on a half-disk has a ``hidden'' symmetry. An attempt to
unveil it is presented in the next section.

\subsection{Proof of Theorem~\ref{thm:main} using a branched double
covering of the disk}\label{subs:proof2}
\begin{figure}[thb!]
\begin{center}
\fbox{\resizebox{0.9\textwidth}{!}{\includegraphics*{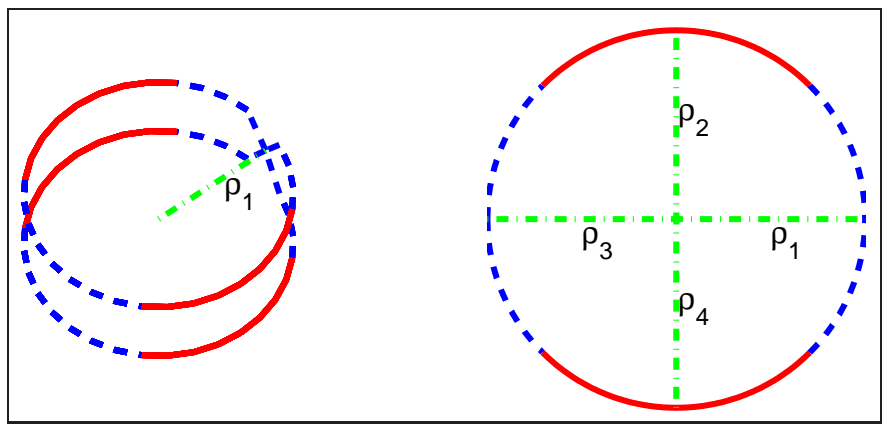}}}
\caption{The problem on the double covering $\D$ (left) of the disk $D$
(right), and the radii on
$D$ which lift to the lines of symmetry for this problem on $\D$.
\label{fig:fig3}}
\end{center}
\end{figure}

Consider an auxiliary spectral problem for  the Laplacian on the
branched double covering $\D$ of the unit disk $D$, with
alternating Dirichlet and Neumann boundary conditions on
quarter-circle  arcs, see Figure~\ref{fig:fig3}. Let
$(r,\phi)$ and $(r,\theta)$ denote  the polar coordinates on
$\D$ and $D$, respectively, with $r\in(0,1]$, $\phi\in[0,4\pi)$, 
$\theta = \phi\pmod{2\pi}\in[0,2\pi)$. 
The metric on $\D\setminus\{O\}$  is a pull-back of the
Euclidean metric from D, where $O=(0,0)$ is the branch point of the
covering. Though $\D$ has a conical singularity at $O$, eigenvalues and eigenfunctions on
$\D$ are well defined using the variational formulation.
Consider  three symmetries of $\D$: $U: (r,\phi) \to (r, 4\pi
-\phi)$, $T: (r, \phi) \to (r, (\phi + 2\pi)\pmod{4\pi})$ and 
$V: (r, \phi) \to (r, (2\pi-\phi)\pmod{4\pi})$. These symmetries are
involutions, they commute with each other, and satisify $V=U\circ T$.
Symmetries $U$ and $V$ are axial symmetries, and $T$ is an
intertwining of sheets of $\D$. By the spectral theorem we find a
basis of eigenfunctions that are either even or odd with respect
to $T$, $U$ and $V$. Consider a space $E_{-}$ of eigenfunctions on
$\D$ that are {\it odd} with respect to $T$ and the corresponding
spectrum $\sigma_{-}(\D)$. We have $E_{-}=E_{-}^{+,-} \cup
E_{-}^{-,+}$, where $E_{-}^{+,-}$ is a subspace of eigenfunctions
that are even with respect to $U$ and odd with respect to $V$, and
$E_{-}^{-,+}$ is a subspace of eigenfunctions that are odd with
respect to $U$ and even with respect to $V$. Denote
$F_U=\{\phi=0\}\cup\{\phi=2\pi\}$ and
$F_V=\{\phi=\pi\}\cup\{\phi=3\pi\}$ the fixed point sets of
$U$ and $V$. Any $f \in  E_{-}^{+,-}$ (respectively, $f\in
E_{-}^{-,+}$) satisfies Neumann  (respectively, Dirichlet)
condition on $F_U$ and Dirichlet (respectively, Neumann)
condition on $F_V$.

Choose a coordinate system on $\D$ in such a way that $\theta=0$
corresponds to the radius $\rho_1 \subset D$ on Figure
\ref{fig:fig3}. For any eigenfunction $f \in  E_{-}^{+,-}$
consider its restriction on the ``upper'' part $\tilde \D = \{
(r,\phi)|\, 0< r < 1, \,\, 0\le \phi < 2\pi\}$. Then $f|_{\tilde
\D}$ projects to an eigenfunction of our boundary problem on a
disk $D$ {\it with a cut} along a diameter $\rho_1\cup \rho_3$: on
$\rho_1$ it satisfies the Neumann condition and on $\rho_3$ it
satisfies the Dirichlet condition. Similarly, any eigenfunction $f
\in  E_{-}^{+,-}$ can be projected from $\tilde \D$ to an
eigenfunction of our boundary problem on a disk $D$ with the same
cut, but now it satisfies Dirichlet condition on $\rho_1$ and
Neumann condition on $\rho_3$. In either case, we obtain an
eigenfunction of the Problem I. Hence, $\sigma_{-}(\D)$ equals
$\sII$ with doubled multiplicities.

Now, let us choose the coordinate system differently so that
$\theta=0$  corresponds to the radius $\rho_2$. Arguing in
exactly the same way as above we obtain that $\sigma_{-}(\D)$
equals $\sII$ with doubled multiplicities. Therefore,
$\sI=\sII$ with account of multiplicities which
completes the proof of the theorem. \qed

\smallskip

\noindent{\bf Remark.} The construction of a ``common'' covering for
Problems I and II described above can be viewed as an application of
Sunada's approach (\cite{Sum}) to mixed Dirichlet-Neumann problems.

\subsection{Proof of Theorem~\ref{thm:main} by Dirichlet-to-Neumann
type mappings}\label{subs:proof3}
\begin{figure}[thb!]
\begin{center}
\fbox{\resizebox{0.4\textwidth}{!}{\includegraphics*{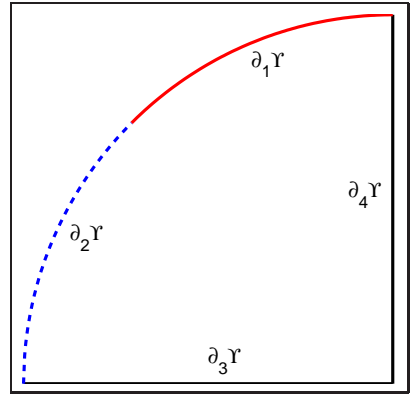}}}
\caption{The quarter-disk $\Upsilon$.
\label{fig:fig2c}}
\end{center}
\end{figure}

For those who prefer operator theory to geometric
constructions we sketch yet another proof of the main theorem.
Consider the following auxiliary problem. Let $\Upsilon$ be a
quarter-disk introduced in
section~\ref{subs:proof1}. Denote $\partial_1\Upsilon=\{|z|=1\,,\ \arg z\in(\pi/2,3\pi/4)\}$,
$\partial_2\Upsilon=\{|z|=1\,,\ \arg z\in(3\pi/4,\pi)\}$, $\partial_3\Upsilon=\{\re z\in(-1,0)\,,\ \im z=0\}$,
$\partial_4\Upsilon=\{\re z=0\,,\ \im z\in(0,1)\}$, so that
$\partial\Upsilon=\overline{\partial_1\Upsilon\cup\partial_2\Upsilon
\cup\partial_3\Upsilon\cup\partial_4\Upsilon}$, see Figure~\ref{fig:fig2c}.

Let, for
a given $\lambda\in\Rbb$, $w(z)$ satisfy the equation
\begin{equation}\label{eq:weq}
-\Delta w = \lambda w\quad\text{in }\tilde\Omega\,,
\end{equation}
and the boundary conditions
\begin{equation}\label{eq:wbc}
w|_{\partial_1\Upsilon}=0\,,\qquad\left.\pd{w}{n}\right|_{\partial_2\Upsilon}=0
\end{equation}
(we do not impose at the moment any boundary conditions on $w$ on
$\partial_{3,4}\Upsilon$). Denote $\xi=w|_{\partial_4\Upsilon}$,
$\eta=\left.\pd{w}{n}\right|_{\partial_4\Upsilon}$,
$p=w|_{\partial_3\Upsilon}$,
$q=\left.\pd{w}{n}\right|_{\partial_3\Upsilon}$, Consider four
linear operators which depend on $\lambda$ as a parameter:
\begin{align}
\DD: \xi \mapsto p, \quad
&\text{subject to} \quad q=0\,,\label{eq:DD}\\
\DN: \xi \mapsto q, \quad
&\text{subject to} \quad p=0\,,\\
\ND: \eta \mapsto p,  \quad
&\text{subject to}\quad q=0\,,\\
\NN:\eta \mapsto q,  \quad
&\text{subject to}\quad p=0\,.\label{eq:NN}
\end{align}
These operators acting on the radius $\partial_4\Upsilon$ are well
defined as long as $\lambda$ does not belong to the spectra of any
of the four homogeneous boundary value problems \eqref{eq:weq},
\eqref{eq:wbc}  with Dirichlet or Neumann boundary conditions
imposed on $\partial_4\Upsilon$ and $\partial_3\Upsilon$ (cf.
\cite{Fried}). Consider an operator
$\mathcal{C}_\lambda=\DD^{-1}\ND\NN^{-1}\DN$. Theorem
\ref{thm:main} then follows from
\begin{prop}\label{pr:DN}
$$
\lambda\in\sI \ \Longleftrightarrow \ \mu=-1
\text{ is  an eigenvalue of } \mathcal{C}_\lambda \ \Longleftrightarrow \ 
\lambda\in\sII\,.
$$
\end{prop}
Proposition \ref{pr:DN} is obtained by re-writing the Problems $\tilde{\text{I}}$ and
$\tilde{\text{II}}$ in terms of operators \eqref{eq:DD}--\eqref{eq:NN}.
We leave the details of the proof to an interested reader. \qed

\section{Extensions, generalizations, open questions}\label{sec:ext}
\subsection{From half-disks to quarter-spheres}\label{subs:qs}
\begin{figure}[thb!]
\begin{center}
\fbox{\parbox{0.4\textwidth}{\resizebox{0.38\textwidth}{!}{\includegraphics*{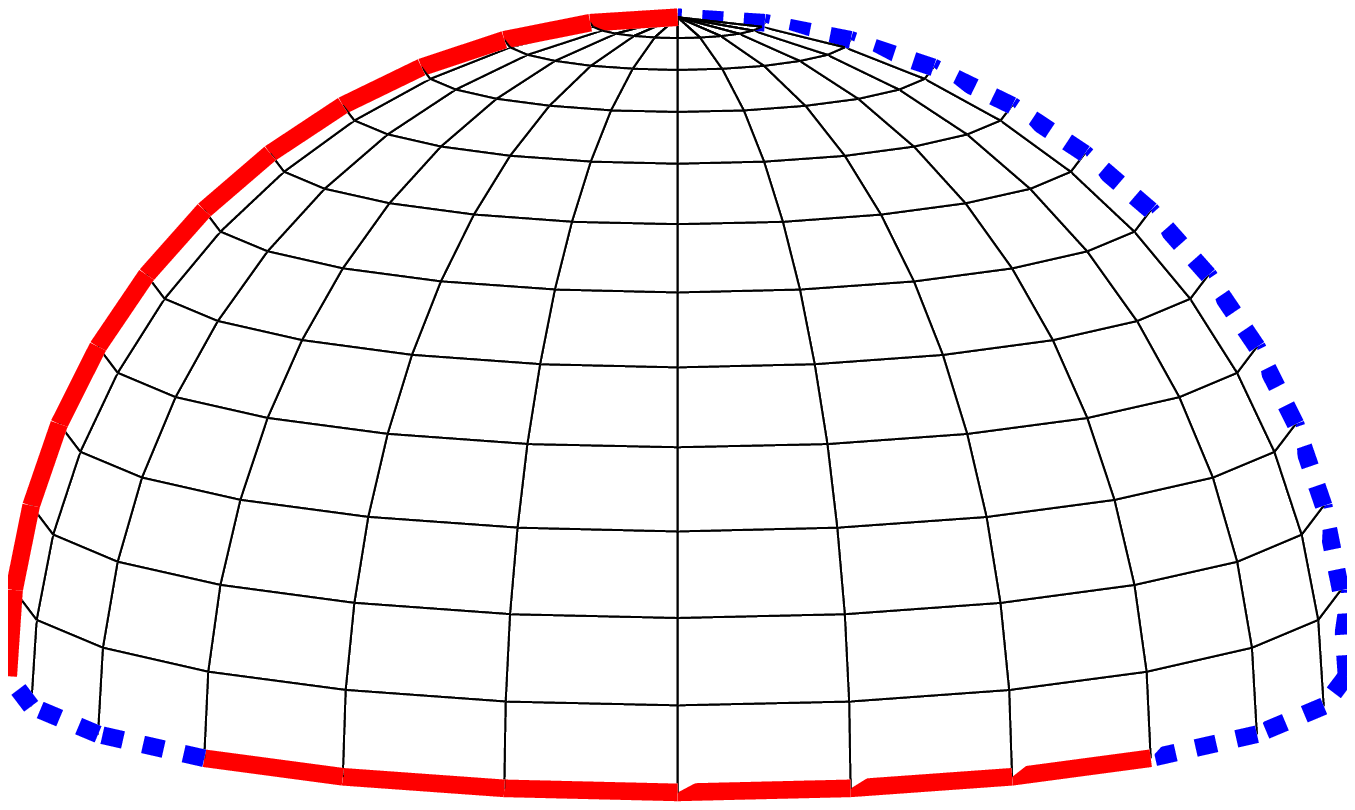}}}\qquad
\parbox{0.4\textwidth}{\resizebox{0.38\textwidth}{!}{\includegraphics*{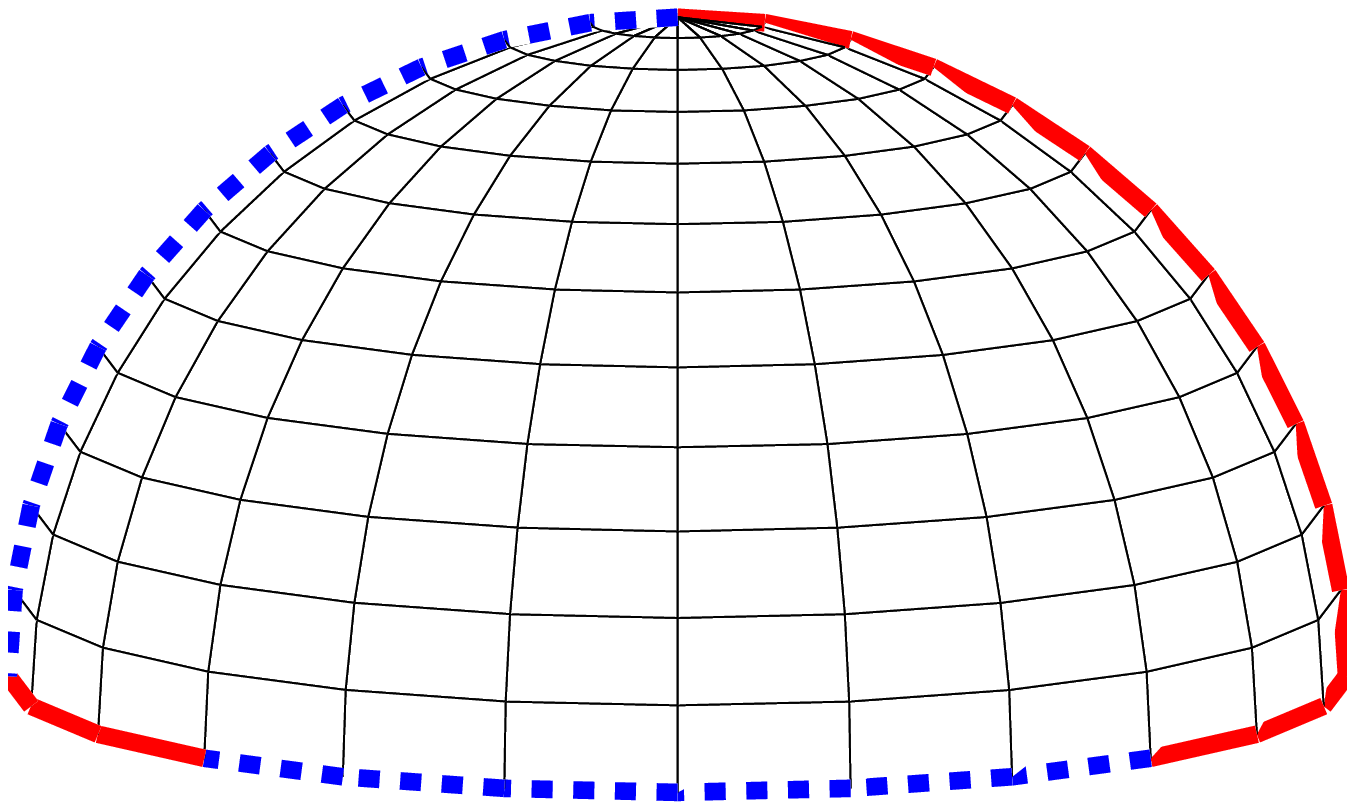}}}}
\caption{Dirichlet-Neumann isospectral problems on quarter-spheres.
The upper semicircles are divided into two equal parts, the lower
semicircles are divided in proportion 1:2:1.
\label{fig:fig6}}
\end{center}
\end{figure}

Consider two quarter-spheres with the boundary conditions as shown
in Figure~\ref{fig:fig6}. 
To prove that they are isospectral one
can use the same trick as shown on Figure~\ref{fig:fig2}. In
general, analogous argument works for half-disks endowed with an
arbitrary radial metric $ds^2=f(|z|)dzd\bar z$ (note that the
matching conditions on Figure~\ref{fig:fig6} are imposed along the
radii), quarter-spheres being a special case for a metric
$\displaystyle  ds^2=\frac{4 dz d\bar z}{(1+|z|^2)^2}$.
Example on Figure~\ref{fig:fig6} was in fact the first 
non-trivial isospectral boundary decomposition that we observed, and it motivated
our study, see section~\ref{subs:g2}.

\subsection{Domains built from sectorial blocks}\label{subs:sb}
Example of section~\ref{subs:ME}  can be also generalized to a class of
domains  constructed by gluing together four copies of a sectorial
block, i.e. a domain bounded by the sides of an acute angle and an
arbitrary continuous curve (without self-intersections) inside it.
Namely, let $0<\alpha<\pi/2$, and choose any points $z_1,z_2\ne 0$
such that $\arg z_1=0$ and $\arg z_2=\alpha$. Now, let $\Gamma_1$
be a piecewise smooth non-self-intersecting curve with end-points
$z_1$, $z_2$ which lies in the sector $\{0<\arg z<\alpha\}$, and
let $K_1$ denote an open set bounded by the radii $[0,z_1]$,
$[0,z_2]$, and the curve $\Gamma_1$.

Let now $\Scal_\beta:(r,\phi)\mapsto (r,2\beta-\phi)$ be a map
which  sends a point into its mirror image with respect to the
axis $\{\arg z=\beta\}$. Let
\begin{equation}\label{eq:Gj}
\Gamma_2:=\Scal_\alpha\Gamma_1\,,\quad \Gamma_3:=\Scal_{2\alpha}\Gamma_2\,,\quad \Gamma_4:=\Scal_{2\alpha}\Gamma_1
\end{equation}
and
$$
K_2:=\Scal_\alpha K_1\,,\quad K_3:=\Scal_{2\alpha}K_2\,,\quad K_4:=\Scal_{2\alpha}K_1\,,
$$
and let $K$ be the interior of $\overline{K_1\cup K_2\cup K_3\cup
K_4}$.  The domain $K$ is bounded by the radii $[0,z_1]$,
$[0,\Scal_{2\alpha}z_2]$ and the curve  $\overline{\Gamma_1\cup
\Gamma_2\cup \Gamma_3\cup \Gamma_4}$

\begin{figure}[thb!]
\begin{center}
\fbox{\resizebox{0.9\textwidth}{!}{\includegraphics*{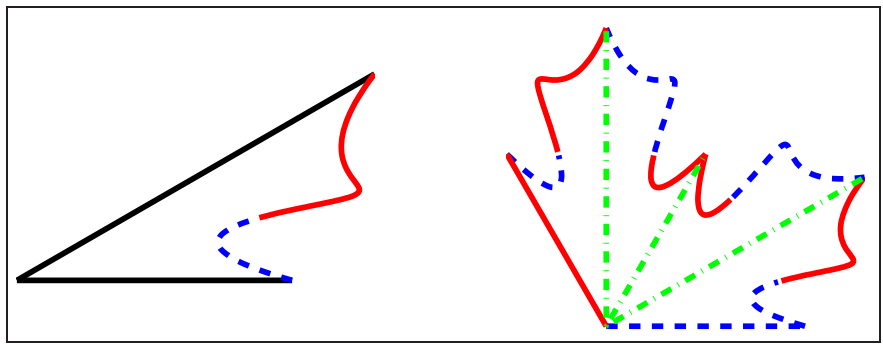}}}
\caption{A sectorial block $K_1$ and the resulting domain $K$. The
spectral problem on $K$ with boundary conditions as shown is
Dirichlet-Neumann isospectral. \label{fig:fig7}}
\end{center}
\end{figure}

We construct a family of pairwise Dirichlet-Neumann isospectral
boundary  value problems on $\Omega$ in the following way. Suppose
$\Gamma_1$ is decomposed into a union of two non-intersecting sets
$\Gamma_{1,1}$ and  $\Gamma_{1,2}$ (one of which may be empty).
We define the sets $\Gamma_{j,m}$, $j=1,2,3,4$; $m=1,2$ similarly
to \eqref{eq:Gj}. We now set
$$
\partial_1 K := \Gamma_{1,1}\cup\Gamma_{2,2}\cup\Gamma_{3,2}\cup\Gamma_{4,1}\cup[0,\Scal_{2\alpha}z_2]
$$
and
$$
\partial_2 K := [0,z_1]\cup \Gamma_{1,2}\cup\Gamma_{2,1}\cup\Gamma_{3,1}\cup\Gamma_{4,2}
$$
(see Figure~\ref{fig:fig7}).

The following result generalises Theorem~\ref{thm:main}.

\begin{thm}\label{thm:gen} With the above notation, the problem
$$
-\Delta u =\lambda u\quad\text{in }K\,,
\quad u|_{\partial_1 K} = 0\,,\quad \left.\pd{u}{n}\right|_{\partial_2 K} = 0\,,
$$
is isospectral with respect to exchanging the Dirichlet and
Neumann  boundary conditions.
\end{thm}

The first proof of Theorem~\ref{thm:main} (see section~\ref{subs:proof1}) is
straightforwardly  adapted for Theorem~\ref{thm:gen}. Note that to
obtain Theorem~\ref{thm:main} we just set $z_1=1$,
$z_2=e^{\ii\pi/4}$, $\Gamma_1=\Gamma_{1,2}=\{e^{\ii t},\
t\in(0,\pi/4)\}$, $\Gamma_{1,2}=\emptyset$ in Theorem~\ref{thm:gen}.

\begin{figure}[thb!]
\begin{center}
\fbox{\resizebox{0.9\textwidth}{!}{\includegraphics*{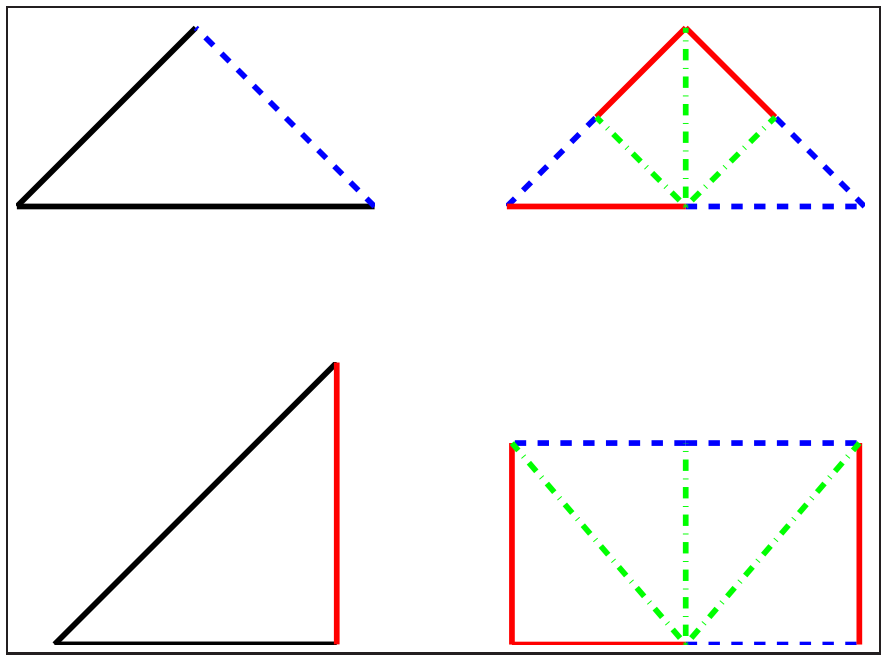}}}
\caption{Two more examples built using sectorial blocks. In the
first  example, $z_1 = 1$, $z_2 = 1/2+\ii/2$,
$\Gamma_1=\Gamma_{1,2}=[z_1,z_2]$, $\Gamma_{1,1}=\emptyset$; the
resulting set $K$ is a triangle. In the second example, $z_1 = 1$,
$z_2 = 1+\ii$,  $\Gamma_1=\Gamma_{1,1}=[z_1,z_2]$,
$\Gamma_{1,2}=\emptyset$; the resulting set $K$ is a $2\times 1$
rectangle. \label{fig:fig8}}
\end{center}
\end{figure}

Other simple examples are illustrated in Figure~\ref{fig:fig8}.

\noindent {\bf Remark.} All our examples of domains admitting
Dirichlet-Neumann isospectrality are constructed using essentially
the same principle. Are there other examples  of such domains? For
instance, all our domains have one axis of symmetry. Do there
exist non-symmetric domains that admit Dirichlet-Neumann
isospectrality? In general, can one characterize in geometric
terms domains admitting Dirichlet-Neumann isospectrality?

\subsection{A necessary condition for Dirichlet-Neumann isospectrality}\label{subs:nes}
After presenting various examples of Dirichlet-Neumann
isospectrality it  would be natural to ask about restrictions.
Intuitively, isospectral decompositions should occur rarely. A
simple necessary condition for a boundary decomposition to be
isospectral is given by
\begin{prop}
\label{heat} If a boundary decomposition $\partial \Omega =
\overline{\partial_1 \Omega \cap \partial_2\Omega}$ of a bounded
planar domain is isospectral with respect to the
Dirichlet-Neumann swap then the total lengths of the parts are
equal: $|\partial_1 \Omega|=|\partial_2\Omega|$.
\end{prop}
\begin{proof} We use asymptotics of the heat trace for a domain with mixed
boundary conditions (see \cite{Gilkey}, \cite{DG}). The first heat
invariant $a_1$ is equal (up to a multiplicative constant) to
$|\partial_N\Omega| - |\partial_D \Omega|$, $\partial_N \Omega$
and $\partial_D \Omega$ being Neumann and Dirichlet parts of the
boundary respectively. This immediately implies the proposition.
\end{proof}

One can probably deduce more sophisticated necessary conditions
for  Dirichlet-Neumann isospectrality using higher heat
invariants.

\subsection{Are there domains not admitting Dirichlet-Neumann isospectrality?}\label{subs:are}
Though in general the question of Kac has a negative answer, there
exist domains that are determined by their Dirichlet spectrum (see
\cite{Zelditch}), for example, a disk. It would be natural to ask
if there are domains not admitting  non-trivial Dirichlet-Neumann
isospectral decompositions of their boundaries.
\begin{conj}
A disk does not admit Dirichlet-Neumann isospectrality.
\end{conj}
We conducted a simple numerical experiment providing some evidence for this conjecture, by
considering boundary decompositions of a unit disk such that
$\partial_1 \Omega$ and $\partial_2 \Omega$ are unions of two
segments each, $|\partial_1\Omega|=|\partial_2 \Omega|=\pi$ by
Proposition \ref{heat}. Each partition is parametrized by a pair
$(\alpha, \beta)$, where $\alpha, \pi - \alpha$ are lengths of
segments in $\partial_1 \Omega$, and $\beta, \pi - \beta$ are
lengths of segments in $\partial_2 \Omega$. For every pair
$(k\pi/24, n\pi/24)$, $0<k\le n<12$ we compute numerically using
FEMLAB \cite{FEM} the $L^2$-norm $\nu(k,n)$ of a vector
$(\lambda_1^\text{I} -\lambda_1^{\text{II}}, \lambda_2^\text{I}
-\lambda_2^\text{II}, \lambda_3^\text{I} -\lambda_3^\text{II})$.
Here $\lambda_i^\text{I}$ are the eigenvalues of the mixed problem
with Dirichlet conditions on $\partial_1 \Omega$ and Neumann
conditions on $\partial_2 \Omega$, and $\lambda_i^\text{II}$ are
the eigenvalues of the problem with the conditions swapped. We
observe that for trivial decompositions ($n=k$) the norm $\nu(k,n)$ is by at least an order of magnitude
smaller than for  any non-trivial decomposition. For example, in a
trivially isospectral case $\nu(12,12)=0.0012$, and in a
non-isospectral case $\nu(11,12)=0.0725$ (this value is in fact
the minimal one achieved among all non-trivial combinations).

\section{Dirichlet-Neumann isospectrality and eigenvalue inequalities}\label{sec:ineq}
\subsection{Genus 2: where did Dirichlet-Neumann isospectrality come from}\label{subs:g2}
In this section we briefly describe our motivation to study
Dirichlet-Neumann isospectrality. It comes, quite unexpectedly,
from a problem to obtain a sharp upper bound for the first
positive eigenvalue $\lambda_1$ of the Laplacian on a surface of
genus $2$. It is known (\cite{YY}, \cite {N1}) that on a surface
$M$ of genus $p$
\begin{equation}\label{upperbd}
\lambda_1\area(M) \le 8\pi \left[\frac{p+3}{2}\right].
\end{equation}
On a surface $\cP$ of genus $2$ this implies
\begin{equation}
\label{sharp}
\lambda_1\area(\cP) \le 16\pi.
\end{equation}
In general (\ref{upperbd}) is not sharp, for example for $\gamma=1$ (\cite{N1}).
In \cite{JLNP} we work towards proving the following
\begin{conj}\label{conj1}
There exists a metric on a surface of genus 2 that attains the
upper  bound in (\ref{sharp}).
\end{conj}

The candidate for the extremal metric is a singular metric of constant
curvature $+1$
that is lifted from a sphere $\Sbb^2$. The surface $\cP$ here is viewed
as a branched double
covering over a sphere with $6$ branching points. The branching points
are chosen to be the intersections of $\Sbb^2$ with the coordinate axes
in ${\Rbb}^3$.  The punctured sphere has an octahedral symmetry group, and
the corresponding hyperelliptic cover corresponds to {\em Bolza's surface}
$w^2=z^5-z$ (known also as the {\em Burnside curve}), and has a symmetry group with $96$ elements (a central
extension by ${\Bbb Z}_2$ of an octahedral group), the largest possible
symmetry group for surfaces of genus $2$, see e.g. \cite{I,KW}.

Note that $\area(\cP)=2\area(\Sbb^2)=8\pi$ and, therefore it suffices to
show that
\begin{equation}
\label{lambda}
\lambda_1 (\cP)=\lambda_1(\Sbb^2)=2.
\end{equation}
It remains to be proved that there exists a first eigenfunction on 
$\cP$ that projects to $\Sbb^2$ , i.e. which is even with respect to 
the hyperelliptic involution $\tau$ intertwining the sheets of the
double cover. We conjecture (see Conjecture~\ref{conj2}) that a first eigenfunction on $\cP$ can not be odd
with respect to $\tau$. The symmetry group of $\cP$ contains many
commuting involutions, and this allows us to exploit the ideas of section~\ref{subs:proof2}. 
Consider an odd eigenfunction (with respect to $\tau$) and symmetrize it with respect to those
involutions. On their fixed point sets we get either Dirichlet or Neumann conditions.
Applying the projection $\cP \to \Sbb^2$ we obtain an eigenfunction of a spectral problem on a sphere
{\it with cuts} along certain arcs of great circles, where Dirichlet or Neumann conditions are imposed.
In particular, in this way we obtain mixed boundary problems shown
on Figure \ref{fig:fig6}. These problems are isospectral, since the
spectrum of each problem coincides with  the odd (with respect to $\tau$) part
of the spectrum of $\cP$. 

All the details of this argument will appear in \cite{JLNP}.

\subsection{Bounds on the first eigenvalue of mixed boundary problems}\label{subs:bounds}
Dirichlet-Neumann isospectrality can be viewed as a special case
of the following question. Consider a mixed Dirichlet-Neumann
problem on a domain with a boundary of length $l$, where the
Dirichlet and the Neumann conditions are specified on parts of the
boundary of total length $l/2$ each. For a given domain, how does
the geometry of the boundary decomposition  affect the spectrum?
We discuss this question in relation to the first eigenvalue $\lambda_1$.

It is natural to ask how large and how small can $\lambda_1$ be.
Extremal boundary decompositions for the first eigenvalue
of a mixed Dirichlet-Neumann problem are studied in
\cite{Denzler1}. In particular, it is  proved that $\lambda_1$ can
get arbirarily close to the first eigenvalue of the pure Dirichlet
problem (which is hence a sharp upper bound for $\lambda_1)$: 
it is achieved in the limit as Dirichlet and Neumann
conditions get uniformly distributed on the boundary. 
It is also shown that a decomposition 
minimizing $\lambda_1$ always exists for bounded Lipschitz domains. 
However, an explicit minimizer is found only for a disk, 
where Dirichlet and Neumann conditions have to be imposed on
half-circles 
\cite{Denzler1}.

The problem of comparing the first eigenvalues for different
boundary decompositions seems to be rather transcendental in
general. Below we present a result,  communicated to us by Brian
Davies and Leonid Parnovski, that applies to a special case of
axisymmetric/centrally symmetric decompositions.

Let $\Xi_a$ be a simply connected planar domain, which is
symmetric with  respect to an axis $d$. We consider a mixed
boundary value spectral problem for the Laplacian on $\Xi_a$ with
some combination of Dirichlet and Neumann boundary conditions on
$\partial\Xi_a$ which is also symmetric with respect to $d$.
Denote the first eigenvalue of this problem by $\lambda_1(\Xi_a)$.

\begin{figure}[thb!]
\begin{center}
\fbox{\resizebox{0.9\textwidth}{!}{\includegraphics*{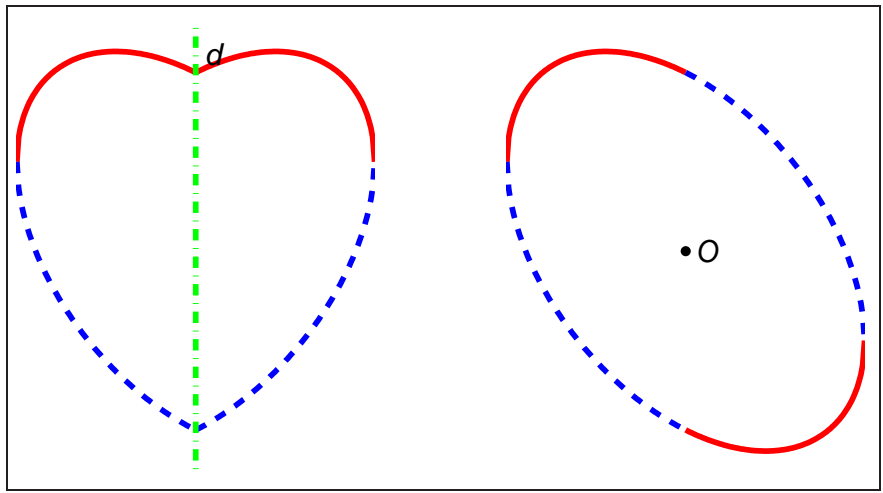}}}
\caption{Axisymmetric domain $\Xi_a$ and centrally symmetric domain $\Xi_c$.
\label{fig:fig9}}
\end{center}
\end{figure}

Let $\Xi_{1,2}$ denote two halves of $\Xi_a$ lying on either side
of  $d$, and let $\tilde{\Xi}_2$ be an image of $\Xi_{1}$ under
the central symmetry with respect to the midpoint $O$ of the
interval $d_\Xi:=\Xi_a\cap d$. Consider a centrally symmetric
domain $\Xi_c$ which is the interior of
$\overline{\Xi_1\cup\tilde{\Xi}_2}$, and the spectral mixed
boundary value problem on $\Xi_c$ with boundary conditions on
$\partial\tilde{\Xi}_2$ centrally symmetric to the ones on
$\partial\Xi_1$, see Figure~\ref{fig:fig9}. Denote the first
eigenvalue of this problem by $\lambda_1(\Xi_c)$.

\begin{thm}\label{thm:DP}{\rm (\cite{DP})} $\lambda_1(\Xi_c)\ge\lambda_1(\Xi_a)$.
\end{thm}

\begin{proof} Consider an auxiliary boundary value problem for the
Laplacain on $\Xi_1$ obtained by keeping the given boundary
conditions on $\partial\Xi_1\setminus d_\Xi$ and imposing the
Neumann condition on $d_\Xi$. Denote the first eigenvalue of this
auxiliary problem by $\lambda_1(\Xi_1)$. By the variational
principle and Dirichlet-Neumann bracketing argument,
$\lambda_1(\Xi_c)\ge \lambda_1(\Xi_1)$. On the other hand, as the
first eigenfunction of the symmetric problem (corresponding to the
eigenvalue $\lambda_1(\Xi_a)$) should be symmetric with respect to
$d$ and therefore should satisfy the Neumann condition on $d_\Xi$,
we have  $\lambda_1(\Xi_a)=\lambda_1(\Xi_1)$, thus implying the
result. Note that the equality $\lambda_1(\Xi_c)=\lambda_1(\Xi_1)$
(and therefore the equality $\lambda_1(\Xi_c)=\lambda_1(\Xi_a)$)
can be attained if and only if $\Xi_1$ has an additional line of
symmetry $d_1$ perpendicular to  $d$ and passing through the
midpoint $O$ of $d_\Xi$, with the boundary conditions being
imposed on $\partial\Xi_1$ symmetrically with respect  to $d_1$.
\end{proof}

\begin{figure}[thb!]
\begin{center}
\fbox{\parbox[t]{0.31\textwidth}{\resizebox{0.3\textwidth}{!}{\includegraphics*{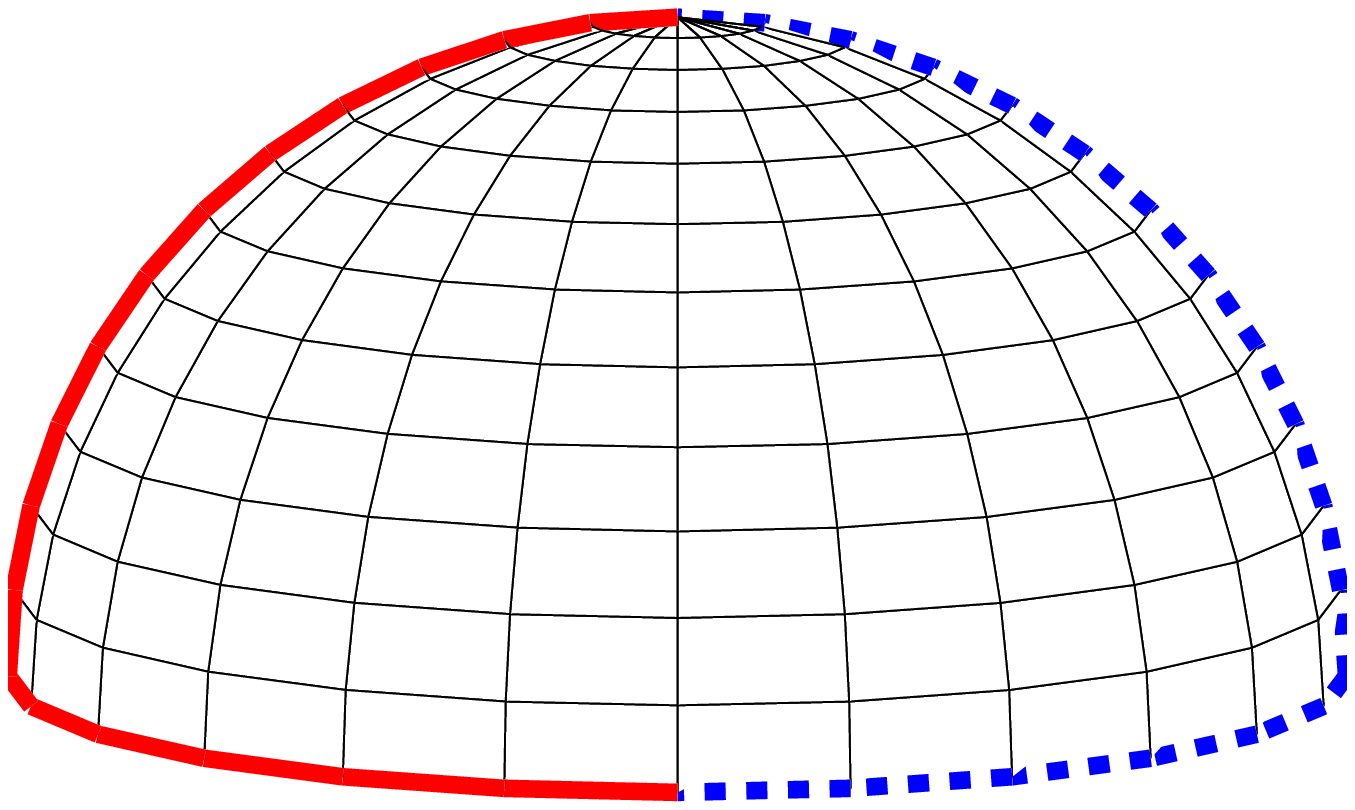}}}\  
\parbox[t]{0.31\textwidth}{\resizebox{0.3\textwidth}{!}{\includegraphics*{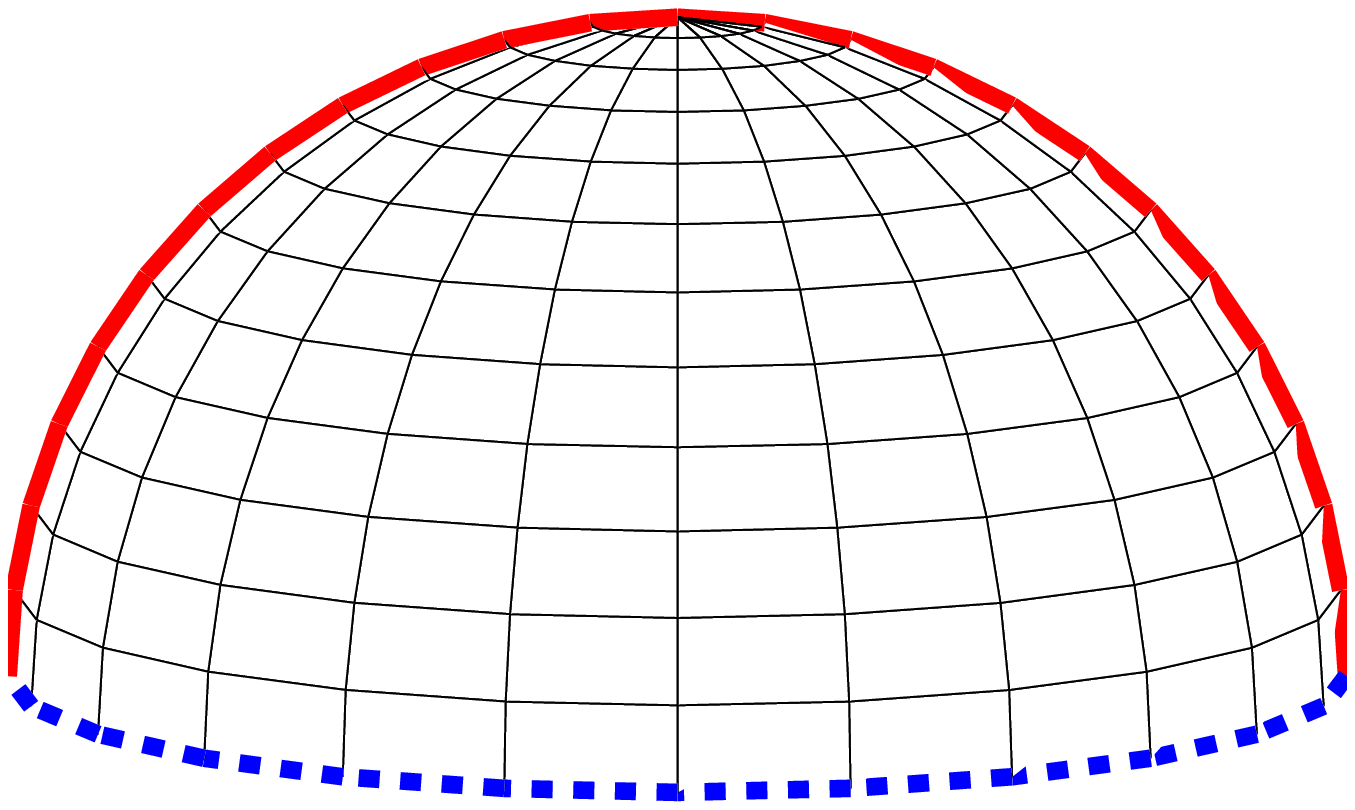}}}\ 
\parbox[t]{0.31\textwidth}{\resizebox{0.3\textwidth}{!}{\includegraphics*{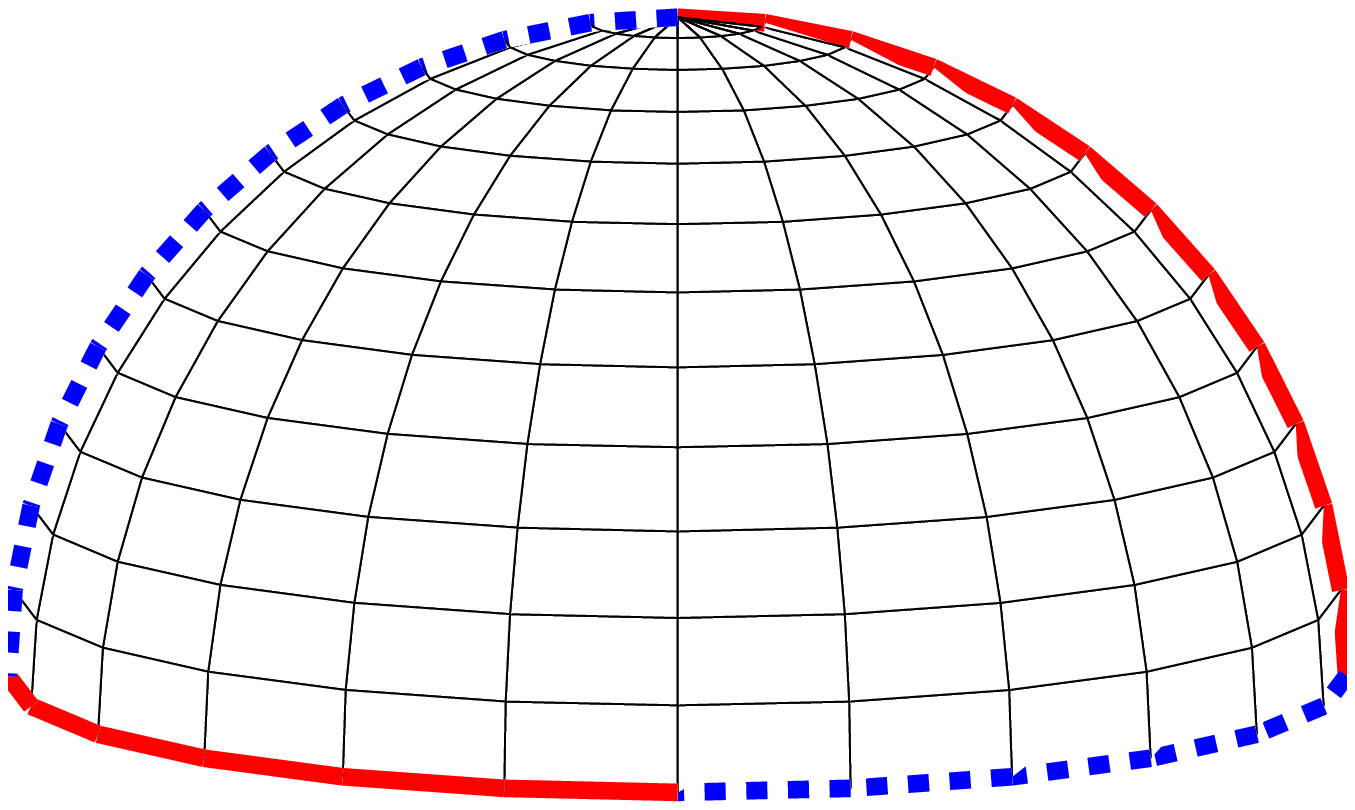}}}}
\caption{Axisymmetric (domains $Q_{a,1}$, left, and $Q_{a,2}$, centre) and centrally symmetric 
(domain $Q_c$, right)
positioning of Dirichlet and Neumann boundary conditions on the
quarter-sphere. The first eigenvalue of the Laplacian is  larger
in the centrally symmetric case: $\lambda_1(Q_c)>\lambda_1(Q_{a,j})$, $j=1,2$.\label{fig:fig12}}
\end{center}
\end{figure}

Theorem \ref{thm:DP} can be used for obtaining estimates of
eigenvalues of boundary value problems on domains
with two lines of symmetry. For example, the boundary of a quarter-sphere has a natural decomposition
into two halves of great circles. If we impose the mixed
Dirichlet-Neumann boundary  conditions on the halves of these
great circles as shown in Figure~\ref{fig:fig12}, we immediately
obtain that the first eigenvalue in any of the axisymmetric cases is
smaller than the first eigenvalue in the centrally symmetric case.

We would like to conclude with another inequality on the first
eigenvalue for quarter-spheres  that we need to check in order to
complete the proof of sharpness of (\ref{sharp}) in \cite{JLNP}.
Let $Q$ be any of the two isospectral problems on a
quarter-sphere shown on Figure \ref{fig:fig6}, and remind that $Q_{a,2}$ is the problem 
shown in the middle of Figure~\ref{fig:fig12} (with the Dirichlet condition imposed 
on one half of the big circle and the Neumann condition on another).
\begin{conj}
\label{conj2}
$\lambda_1(Q)>\lambda_1(Q_{a,2})$
\end{conj}
One can immediately check that $\lambda_1(Q_{a,2})=\lambda_1(\Sbb^2)=2$.
An affirmative solution of Conjecture \ref{conj2} excludes the
possibility that the first eigenfunction on $\cP$ is odd with
respect to the intertwining  of sheets (see section~\ref{subs:g2}), and
hence we have
\begin{thm}{\rm (\cite{JLNP})} Conjecture \ref{conj2} implies
Conjecture \ref{conj1}.
\end{thm}
Using FEMLAB \cite{FEM} one can verify Conjecture \ref{conj2}
numerically: $\lambda_1(Q) \approx 2.28 > 2$. Our current
project is to find a rigorous (possibly, computer-assisted) proof
of this conjecture.

\medskip

\noindent {\it Acknowledgements.} We are grateful to Leonid Parnovski for
many useful suggestions, in particular those providing the proof
in section~\ref{subs:proof2}. We also thank
Leonid Polterovich for helpful remarks on a preliminary version of this paper.

\end{document}